\documentclass[10pt]{amsart}
\newcommand{\vs}[1]{\vspace{#1}}

\newcommand{\ra}{\rightarrow}
\newcommand{\R}{\mathbb{R}}

\newcommand{\N}{\mathbb{N}}

\newcommand{\ep}{\epsilon}

\newcommand{\8}{\infty}
\newcommand{\nn}{\nonumber}
\newcommand{\ee}{\end{eqnarray}}
\newcommand{\be}{\begin{eqnarray}}
\newcommand{\dfr}{\mbox{{\rm Diff}$_\mu^r(M)$}}

\newtheorem{thm}{Theorem}[section]
\newtheorem{lem}[thm]{Lemma}
\newtheorem{prop}[thm]{Proposition}

\begin{document}

\title{Hyperbolic Invariant Sets With Positive Measures}
\author{Zhihong Xia}
\address{Department of Mathematics \\ Northwestern University \\
Evanston, Illinois 60208}
\thanks{Research supported in part by National Science Foundation.}
\date{April 10, 2005, revised version}
\email{xia@math.northwestern.edu}

\begin{abstract}
In this short paper we prove some results concerning
volume-preserving Anosov diffeomorphisms on compact manifolds. The
first theorem is that if a $C^{1 + \alpha}$, $\alpha >0$,
volume-preserving diffeomorphism on a compact connected manifold
has a hyperbolic invariant set with positive volume, then the map
is Anosov. This is not necessarily true for $C^1$ maps. The proof
uses a Pugh-Shub type of dynamically defined measure density
points, which are different from the standard Lebesgue density
points. We then give a {\it direct} proof of the ergodicity of
$C^{1+\alpha}$ volume preserving Anosov diffeomorphisms, without
using the usual Hopf arguments or the Birkhoff ergodic theorem.
The method we introduced also has interesting applications to
partially hyperbolic and non-uniformly hyperbolic systems.
\end{abstract}

\maketitle

\section{Introduction and statement of main results}

We consider volume-preserving or symplectic diffeomorphisms on a
compact connected Riemannian manifold $M$. Let $\dfr$ be the the set
of all $C^r$ diffeomorphisms preserving a smooth volume $\mu$ on
$M$. If $r$ is not an integer, $r= k + \alpha$ for some positive
integer $k$ and $0 < \alpha < 1$, it is understood that the functions
in $\dfr$ are $C^k$ functions with $\alpha$-H\"older $k$-th
derivatives.

An invariant set $\Lambda \subset M$ is said to be {\em hyperbolic}\/
if there is a continuous splitting of $T_xM = E^s_x \oplus E^u_x$ for
every $x \in \Lambda$ and constants $C>0$, $\lambda >1$ such that
\be
df_x(E^s_x) &=& E^s_{f(x)} \; \mbox{ and } \; df_x(E^u_x) = E^u_{f(x)}
\nn \\ |df^n_x v^s_x| &\leq& C \lambda^{-n} |v^s|,\; \mbox{ for all } v^s
\in E^s_x, \; n \in \N \nn \\ |df^{-n}_x v^u_x| &\leq& C \lambda^{-n}
|v^u|,\; \mbox{ for all } v^u \in E^u_x, \; n \in \N \nn \ee

If the whole manifold $M$ is hyperbolic for some $f \in \dfr$, then
$f$ is said to be {\em Anosov}. Not all manifolds can support
Anosov diffeomorphisms.

Typical examples of hyperbolic invariant sets are Cantor sets as in
Smale's horseshoe map. The following simple proposition explains why
this is the case.

\begin{prop}
Let $f \in \dfr$, $r \geq 1$, be a volume-preserving
diffeomorphism on a compact manifold $M$. Let $\Lambda \subset M$ be a
closed hyperbolic invariant set. If the interior of $\Lambda$ is
non-empty, then $f$ is Anosov on $M$ and $\Lambda =M$.
\label{prop}
\end{prop}

This proposition and its simple proof, given in the next section, will
motivate our main result of this paper. The proof uses the fact that
the recurrent points are dense on the manifold. This is a consequence
of the volume-preserving property.  Without the volume-preserving or
the dense recurrent points condition, the proposition is not true, we
refer to Fisher \cite{Fisher04} for a counter-example. Fisher also
give a proof of the above proposition. On the other hand, it is an
open problem whether there are any Anosov diffeomporphisms with
wandering domains.

A natural question one asks is whether there is any hyperbolic
invariant set with a positive measure for a volume-preserving
non-Anosov diffeomorphism. The answer is yes for $C^1$
diffeomorphisms, as Bowen's example of fat horseshoe shows
\cite{Bowen75}, see also Robinson \& Young \cite{RY80}. However, if
the map is assumed to be $C^{1 + \alpha}$ for some $\alpha >0$, then
the answer is no. This is the main result of this paper.

\begin{thm}
Let $f \in \dfr$, $r>1$, be a volume preserving diffeomorphism on a
compact manifold $M$. Let $\Lambda \subset M$ be a closed hyperbolic
invariant set. If $\mu(\Lambda) >0$, then $f$ is Anosov on $M$ and
$\Lambda = M$.
\label{thm}
\end{thm}

It is not surprising that the map is required to be $C^{1+
  \alpha}$. As various examples show, measure-theoretical properties
  are often not respected by $C^1$ maps. Additional smoothness, even
  though very little, guarantees certain regularities in measure.

Our proof uses a special type of measure density points different from
the Lebesgue density points. The density basis for our density points
are dynamically defined. It is similar to the juliennes defined by
Pugh \& Shub \cite{PS00} \cite{PS03}. But our case is much simpler.

Our method also provides a direct proof of the ergodicity of $C^{1 +
  \alpha}$ volume-preserving Anosov diffeomorphisms, without using the
Hopf arguments or the Birkhoff ergodic theorem. However, we do use the
absolute continuity of stable and unstable foliations. This is given
in the last section of the paper.

Another result of this paper Lemma \ref{lemf} is of interest in its
own right. We showed that for a $C^{1+ \alpha}$ hyperbolic or
partially hyperbolic volume-preserving diffeomorphism, if a set is
invariant, then almost every point of the stable manifold (or unstable
manifold) of almost every point is in the invariant set. This is also
true for non-unifomly hyperbolic invariant set in Pesin theory. This
result can find its applications in various other problems.

\vspace{1ex} We are grateful to M. Viana for pointing out that
Theorem \ref{thm} was also proved, with a different method, by
Bochi and Viana \cite{BV03}. We are also grateful to F. Ledrappier
for showing us another possible proof of Theorem \ref{thm}.

\section{Proof of the proposition}

In this section, we give a simple proof of the Proposition \ref{prop}.

Let $U$ be the interior of the hyperbolic invariant set $\Lambda$. By
the assumption of the proposition, $U \neq \emptyset$. Clearly, $U$ is
invariant. We want to prove that the closure of $U$, $\bar{U}$ is the
whole manifold. We know that $\bar{U}$ is closed, it suffices to show
that $\bar{U}$ is also open.

For any $x \in \bar{U}$, there exists a sequence of points $x_n \in
U$, $n \in \N$, such that $x_n \ra x$ as $n \ra \8$. As $f$ is volume
preserving, by Poincar\'e recurrence theorem, almost every point is
both forward and backward recurrent. Moreover, the set of periodic
points is dense in $U$, since $U$ is hyperbolic. We may choose
$\{x_n\}_{n \in \N}$ to be periodic points. Since $U$ is invariant and
each $x_n$ is an interior point in $U$, then $W^s(x_n)$ and $W^u(x_n)$
are in $U$ for all $n \in \N$. For each fixed $\delta >0$ small, let
$W^s_\delta(x)$ and $W^u_\delta(x)$ be, respectively, the local stable
manifold and unstable manifold of $x$. As $x_n \ra x$, as $n \ra \8$,
we have that $W^s_\delta(x_n) \ra W^s_\delta(x)$ and $W^u_\delta(x_n)
\ra W^u_\delta(x)$ as $n \ra \8$.  This implies that each point on
$W^s_\delta(x)$ or on $W^u_\delta(x)$ is also in the closure of
$U$. Let $y \in W^s_\delta(x)$ and $z \in W^u_\delta(x)$, the same
argument shows that $W^u_\delta(y)$ and $W^s_\delta(z)$ are both in
the closure of $U$. Consequently,
$$W^u_\delta(y) \cap W^s_\delta(z) \in \bar{U},$$ i.e., $\bar{U}$ has
the product structure. This implies that $x$ is in the interior of
$\bar{U}$. Consequently, the set $\bar{U}$ is open. Since $\bar{U}$ is
also closed and $M$ is connected, we have $\bar{U} = M$. i.e., $f$ is
hyperbolic on $M$.

This proves the proposition.

\section{Proof of the Theorem}

The proof of Theorem \ref{thm} uses a similar idea to the proof of
Proposition \ref{prop}, but the details are much more
complicated. Here the interior points are replaced by density
points. One may regard the density points as measure theoretical
interior points for a set with positive measure.

We need some preliminary results from standard smooth ergodic theory.
It is well-known that the stable and unstable foliations for a $C^1$
Anosov diffeomorphism may not be absolutely continuous. However, for
$C^{1 +\alpha}$ diffeomorphisms, these foliations are absolutely
continuous (Anosov \cite{Anosov67}). Moreover, the stable and unstable
foliations over a hyperbolic (even non-uniformly, cf Pesin
\cite{Pesin77}) invariant set are also absolutely continuous for $C^{1
  + \alpha}$ diffeomorphisms. In fact, the absolute continuity of the
foliations is proved by showing that the holonomy maps of these
foliations are absolutely continuous.

We also need some results on density basis and density points of a
measurable set. Let $A \subset \R^n$ be a measurable set with the
standard Lebesgue measure $m$. A point $x \in \R^n$ is said to be a
Lebesgue density point if $$\lim_{\ep \ra 0} \frac{m(B(x, \ep) \cap
  A)}{m(B(x, \ep))} =1,$$
where $B(x, \ep)$ is the $\ep$-ball
centered at $x$. Lebesgue density theorem states that almost every
point of $A$ is a density point for $A$.

To prove our theorem, we need a different definition of density point.
The Lebesgue density point is defined by a basis of $\ep$-balls. We
replace it by a dynamically defined basis. Let $\Lambda \subset M$ be
a hyperbolic invariant set, we first define a basis on the unstable
manifold for each point $\Lambda$.

For a fixed small real number $\delta >0$, let $W^u_\delta(x)$ be the
local unstable manifold of a point $x \in \Lambda$.  Let $\mu_u$ and
$\mu_s$ respectively be the induced measures of the smooth volume form
$\mu$ on the unstable leaves and stable leaves. Let $n_u$ and $n_s$
respectively be the dimensions be the unstable and stable leaves. For
any positive integer $k$, let $B^u_k(x)$ be a subset of $W^u(x)$
defined by
$$B^u_k(x) = f^{-k}(W^u_\delta(f^k(x))). $$
Clearly, the cubes
$B^u_k(x)$, $k \in \N$ shrinks to the point $x$ as $k \ra \8$. We call
the collection of the sets $\{B^u_k(x) \; | \;  k \in \N, \; x \in
\Lambda \}$ the unstable density basis.

Similarly, we can define the stable density basis $\{B^s_k(x) \; | \;
k \in \N, \; x \in \Lambda \}$, by defining
$$B^s_k(x) = f^{k}(W^u_\delta(f^{-k}(x))). $$ The density basis we
defined has infinite eccentricity.

A point $x \in \Lambda$ is said to be a dynamical density point, or
simply density point, on the unstable foliation if
$$ \lim_{k \ra \8} \frac{\mu_u(B^u_k(x) \cap
\Lambda)}{\mu_u(B^u_k(x))} =1.$$ Similarly, we can define the
dynamical density points on the stable foliation.

\begin{prop}
  The set of points on $\Lambda$ that are both density points on the
  stable foliations and unstable foliations has the full measure in
  $\Lambda$.
\label{propd}
\end{prop}

Our definitions of density points can be regarded as simplified
versions of the juliennes density points defined by Pugh \& Shub
\cite{PS00} \cite{PS03}. The Pugh-Shub density points are defined by
the julienne density basis and they are much more complicated than
what we have here. The above proposition follows from the proof for
the Pugh-Shub's juliennes density point. The proof itself is
similar to the proof of the Lebesgue Density Theorem. The key
properties for the density basis are scaling and engulfing defined as
follows.

(a) Scaling: for any fixed $k \geq 0$, $m(B^u_n(x)) /m( B^u_{n+k}(x))$
    is unformly bounded as $n \ra \8$.

(b) Engulfing: there is a unifom $L$ such that $$B^u_{n+L}(x) \cap
    B^u_{n+L}(y) \neq \emptyset \; \Rightarrow \; B^u_{n+L}(x) \cup
    B^u_{n+L}(y) \subset B^u_n(x). $$

These two properties can be easily verified.

We remark that our density points are defined on the stable and
unstable foliations, we freely used the fact that the stable and
unstable foliations are absolutely continuous.

Let $A$ be a subset of $\Lambda$ such that for any $x \in A$, $x$ is a
density point of $\Lambda$ on both stable foliation and unstable
foliation; and $x$ is a recurrent point, both forward and backward. By
Poincar\'e recurrence theorem, Proposition \ref{propd} and the
absolute continuity of the foliations, the set $A$ has the full
measure in $\Lambda$.

The following is the main lemma in proving our theorem.

\begin{lem}
  Assume that $f$ is a $C^{1+\alpha}$ volume-preserving
  diffeomorphism, for some positive number $\alpha >0$. Fix $x \in A$
  and a positive number $\delta >0$. Then for any $\ep >0$, there
  exists a positive integer $k_0$, depending on $x$ and $\ep$, such
  that for $k \geq k_0$, $$\mu_s(W^s_\delta(f^{-k}(x)) \cap A) \geq
  (1-\ep) \mu_s(W^s_\delta(f^{-k}(x)))$$
  and
  $$\mu_u(W^u_\delta(f^k(x)) \cap A) \geq (1-\ep)
  \mu_u(W^u_\delta(f^k(x)))$$
  i.e., for sufficiently large $k$, the
  set $A$ has a very high density in $W^s_\delta(f^{-k}(x))$ and
  $W^u_\delta(f^k(x))$.
\end{lem}

\noindent {\it Proof of the lemma}: We first prove the lemma for the
unstable foliation. Let $n_u$ be the dimension of the leaves of the
foliation, local unstable manifold $W^u_\delta(x)$ can be identified
with a cube in $E^u_x = \R^{n_u}$ by the exponential map from $E^u_x$
to $W^u(x)$.  Since the leaves of the unstable foliation is smooth,
the conditional measure $\mu_u$ are smoothly equivalent to the
standard Lebesgue measure $m$ on $\R^{n_u}$. i.e., for any point $x
\in \Lambda$, there is a smooth function $g_u(y)$ defined for $y \in
W^u_\delta(x)$ on the local unstable manifold, uniformly bounded away
from zero and infinity, such that
$$\mu_u(E) = \int_{E} g_u dm,$$
where $E$ is a measurable set in
$W^u_\delta(x)$ and $m$ is the standard Lebesgue measure in
$\R^{n_u}$.

For any positive integer $k$, we want to estimate the measure of the
set $W^u_\delta(f^k(x)) \cap A$. Let $B^k_0 =
f^{-k}(W^u_\delta(f^k(x)))$, obviously $B^k_0 \subset W^u_\delta(x)$.
In fact, $B^k_0$ is the set $B^u_k(x)$ in our definition of density
basis. We iterate $B^k_0$ under $f$ and obtain a sequence of sets
$B^k_i = f^i(B^k_0)$, for $i=1, 2, \ldots, k$. The last set in the
sequence is $B^k_k = W^u_\delta(f^k(x))$.

Let $\eta_k = 1- \mu_u(B^k_0 \cap A) / \mu_u(B^k_0)$, then $0 \le \eta_k
\leq 1$. As $x$ is a density point on the unstable foliation,
$$\lim_{k \ra \8} \frac{\mu_u(B^k_0 \cap A)}{\mu_u(B^k_0)} =1,$$
The
number $\eta_k$ is small for large $k$, and $\lim_{k \ra \8} \eta_k
=0$.  Since $f$ is $C^{1 + \alpha}$, there exists a constant $C_1 > 0$
such that $||df_y - df_z|| \leq C_1 |y-z|^\alpha$. Here we abuse the
notation a little by writing $|y-z|$ as the distance between $y$ and
$z$.

Let $\rho^k_i$ be the maximum distance from $f^i(x)$ to the boundary
of $B^k_i$, i.e., $$\rho^k_i = \max_{y \in B^k_i} \{ d(f^i(x), y) \}.
$$
For any $y \in B^k_0$, $||df_y - df_x|| \leq C_1(\rho^k_0)^\alpha$.
In general, for any $y \in B^k_i$, $||df_y - df_{f^{i}(x)}|| \leq C_1
(\rho^k_i)^\alpha$.

Let $J_u(y)=|\det(df_y|_{E^u_y})|$ be the Jacobian of the map $f$ at $y$
restricted on the unstable manifold of $y$. Then $|J_u(x)-J_u(y)| \leq
C_2 |x-y|^\alpha$, for some positive constant $C_2 >0$.

Let $D_0= B^k_0 \backslash A$ and $D_i = B^k_i \backslash A$, for $i
=1, 2, \ldots, k$. These are the complements of $A$ in $B^k_i$. By the
definition of $\eta_k$, $\mu_u(D_0) = \eta_k \mu_u(B^k_0)$. We need to
estimate the measure of $D_1$. For any set $E \subset B^k_0$,
$$\mu_u(f(E)) = \int_{f(E)} g_u dm =\int_{E} J_u (g_u\cdot f) dm =
\int_{E} J_u (g_u\cdot f) g_u^{-1} d\mu_u$$

Since the functions $g_u$ and $g_u^{-1}$ are smooth on any
unstable manifold, the integrand
in the above integral is $C^{1+\alpha}$, there is a constant $C_3 >0$
such that $$|J_u(y) (g_u\cdot f)(y) g_u^{-1}(y) - J_u(x) (g_u\cdot f)(x)
g_u^{-1}(x)| \leq C_3 |x - y|^\alpha , $$ for all $x\in \Lambda$, $y \in
W^u_\delta(x)$. Therefore,
$$ |\mu_u(f(E)) - (J_u(x) (g_u\cdot f)(x)
g_u^{-1}(x)) \mu_u(E)| \leq C_3 (\rho^k_0)^\alpha \mu_u(E).$$
Consequently,
$$\mu_u(D_1) \leq (J_u(x) (g_u\cdot f)(x)
g_u^{-1}(x)) \mu_u(D_0) + C_3 (\rho^k_0)^\alpha \mu_u(D_0)$$
and $$\mu_u(B^k_1\backslash D_1) \geq (J_u(x) (g_u\cdot f)(x)
g_u^{-1}(x)) \mu_u(B^k_0 \backslash D_0) - C_3 (\rho^k_0)^\alpha
\mu_u(B^k_0\backslash D_0)$$
and therefore
$$\mu_u(D_1) \leq \eta_k \frac{(1 + C_3
(\rho^k_0)^\alpha)}{(1 - C_3
(\rho^k_0)^\alpha)} \mu_u(B^k_1).$$ By induction on $i$, we have
$$\mu_u(D_k) \leq \eta_k \frac{(1 + C_3
(\rho^k_0)^\alpha)}{(1 - C_3
(\rho^k_0)^\alpha)}\frac{(1 + C_3
(\rho^k_1)^\alpha)}{(1 - C_3
(\rho^k_1)^\alpha)} \cdots \frac{(1 + C_3
(\rho^k_{k-1})^\alpha)}{(1 - C_3
(\rho^k_{k-1})^\alpha)}\mu_u(B^k_k). $$

The map $df: T_{\Lambda}M \ra T_{\Lambda}M$ uniformly expands vectors
on the unstable splitting. That uniform expansion extends to local
unstable manifolds $W^u_\delta(x)$, $x \in \Lambda$ if $\delta$ is
chosen small enough. There exist positive real numbers $C_4 > C >0$
and $\lambda> \lambda_1 >1$ (here $C$ and $\lambda$ are the same as
those in the definition of the hyperbolic invariant set) such that
$\rho^k_{k-1} \leq C_4 \lambda_1^{-1} \delta$ and $\rho^k_{i} \leq C_4
\lambda_1^{-(k-i)} \delta$, for $i=0, 1, \ldots, k-1$. This implies
that \be \mu_u(D_k) & \leq & \eta_k \frac{(1 + C_3
(C_4 \lambda_1^{-k}\delta)^\alpha)}{(1 - C_3
(C_4 \lambda_1^{-k}\delta)^\alpha)}\frac{(1 + C_3
(C_4 \lambda_1^{-k+1}\delta)^\alpha)}{(1 - C_3
(C_4 \lambda_1^{-k+1}\delta)^\alpha)} \nn \\
&& \cdots \frac{(1 + C_3
(C_4 \lambda_1^{-1}\delta)^\alpha)}{(1 - C_3
(C_4 \lambda_1^{-1}\delta)^\alpha)}\mu_u(B^k_k) \nn \\
 &<& \eta_k \mu_u(B^k_k) \prod_{i=1}^\8 (1+
C_3(C_4\lambda_1^{-i}\delta)^\alpha )/ \prod_{i=1}^\8
(1- C_3(C_4\lambda_1^{-i}\delta)^\alpha ) \nn\ee
Since $\lambda > 1$, then $\lambda^\alpha >1$, the infinite products
converge.  We have
$$\mu_u(D_k) < C_5 \eta_k \mu_u(B^k_k).$$ for some constant $C_5 >0$.

Choose a positive integer $k_0$ such that for $k \geq k_0$, $\eta_k <
\ep /C_5$, then we have
$$\mu_u(W^u_\delta(f^k(x))
\cap A) \geq (1-\ep) \mu_u(W^u_\delta(f^k(x))),$$
for all $k \geq k_0$.

This proves the statement of the lemma on the unstable foliation. The
part on the stable foliation can be proved in the same way by
considering $f^{-1}$.

This proves the lemma. \vs{1ex}

We return to the proof of the theorem. Let $E$ be the closure of $A$
in $M$. We claim that if $y \in E$, then $W^s_\delta(y) \subset E$ and
$W^u_\delta(y) \subset E$.

Suppose that this is not true. i.e., there is a point $z \in W^u(y)$
and a small $\ep_1$-ball around $z$, $B(z, \ep_1) \subset M$ such that
$B(z, \ep_1) \cap A = \emptyset$. Consequently, there are constants
$\ep_2>0$, depending on $\ep_1$, and $C_6 >0$, independent of $\ep_1$,
such that if $x \in \Lambda$, $|x-y| \leq \ep_2$, then
$$\mu_u(W^u_\delta(f^k(x)) \cap A) < (1- C_6 \ep^{n_u}_1)
\mu_u(W^u_\delta(f^k(x))),$$ for all $k \in \N$.

On the other hand, since $y \in E$, there is a sequence of points $x_i
\in A$, $i \in \N$ such that $x_i \ra y$ and $i \ra \8$. For the above
$\ep_1$, there is positive integer $i_0$ such that if $i \geq i_0$,
$d(x_i, y) \leq \ep_1/3$. For any fixed $\ep >0$, by the lemma above
and the recurrence of $x_i$, there is a positive integer $k >0$ such
that $d(x_i, f^k(x_i)) \leq \ep_1 /3$ and
$$\mu_u(W^u_\delta(f^k(x_i)) \cap A) \geq (1-\ep)
\mu_u(W^u_\delta(f^k(x_i))).$$ Since $d(y, f^k(x_i)) \leq
\frac{2\ep_2}{3} < \ep_2$, choosing $\ep = C_6 \ep^{n_u}_1$ leads to a
contradiction. This contradiction show that if $y \in E$, then
$W^u_\delta(y) \subset E$. Similarly by considering $f^{-1}$, we have
$W^s_\delta(y) \subset E$.

To conclude our proof, for any $y \in E$, $W^u_\delta(y) \subset E$
and $$V = \bigcup_{z \in W^u_\delta(y)} W^s_\delta(z) \subset E.$$
Since $V$ is hyperbolic, $y$ is in the interior of $V$. This implies
that $E$ is an open set. But $V$ is also closed and non-empty. The
connectness of $M$ implies that $E=M$. This implies that $f$ is
Anosov.

Finally, the reason that $\Lambda=M$ in the first place is that $f$ is
ergodic, any invariant set with positive measure must have full
measure and its closure must be the whole manifold.

This proves the theorem.

\section{Ergodicity of volume preserving Anosov diffeomorphisms}

In this section, we give a direct proof of the ergodicity of $C^{1 +
  \alpha}$ volume-preserving Anosov diffeomorphisms, without using the
  usual Hopf arguments or the Birkhoff ergodic theorem.

We first prove the following lemma.

\begin{lem}
Let $f \in \dfr$, $r>1$, be a volume preserving diffeomorphism on a
compact manifold $M$. Let $\Lambda \subset M$ be a
hyperbolic invariant set (not neccessarily closed). If $\mu(\Lambda)
>0$, then for a.e. $x \in \Lambda$, $W^s(x) \subset \Lambda$ modulus a
$\mu_s$-measure zero set in $W^s(x)$ and $W^u(x) \subset \Lambda$ modulus a
$\mu_u$-measure zero set in $W^u(x)$, where $\mu_s$ and $\mu_u$ are
respectively the induced measures of $\mu$ on the stable and unstable
foliations.
\label{lemf}
\end{lem}

\noindent {\it Proof}: For any fixed positive integer $k$
and positive number $\eta>0$, let $\Lambda_{(\eta, k)} \subset \Lambda$
be the set such that $$\frac{\mu_u(B^u_i(x) \cap
\Lambda)}{\mu_u(B^u_i(x))}  > (1 -\eta), \; \mbox{ for all } i \geq k$$
Since almost every point of $\Lambda$ is a density point of $\Lambda$,
for any $\eta>0$, $$\lim_{k\ra \8} \mu(\Lambda_{(\eta, k)}) =
\mu(\Lambda).$$

By Poincar\'e recurrence theorem, for a.e. $x \in \Lambda_{(\eta,
k)}$, there exists a sequence of integers $n_i \ra \8$ such that
$f^{-n_i}(x) \in \Lambda_{(\eta, k)}$. By the distortion estimates
from the last section, $$\mu_u(W^u_\delta(x) \cap \Lambda) \geq
(1-C_5\eta) \mu_u(W^u_\delta(x)),$$ where $C_5$ and $\delta$ are fixed
constants from the last section. Since the above estimate is independent
of $k$, it must hold for almost all $x \in \cup_{k=1}^\8 \Lambda_{(\eta,
k)}$. Since $\mu(\cup_{k=1}^\8 \Lambda_{(\eta, k)}) = \mu(\Lambda)$
for any fixed $\eta >0$, this implies that for a.e. $x \in \Lambda$,
$$\mu_u(W^u_\delta(x) \cap \Lambda) \geq
(1-C_5\eta) \mu_u(W^u_\delta(x)).$$
This is true for all $\eta >0$, therefore, for a.e. $x \in \Lambda$,
$$\mu_u(W^u_\delta(x) \cap \Lambda) = \mu_u(W^u_\delta(x)).$$

This proves the lemma for the unstable foliations. The results on the
stable foliation can be proved in the same way by considering
$f^{-1}$.

This proves the lemma.

\vs{1ex}
Now we can prove the following ergodicity theorem.

\begin{thm}
Let $f \in \dfr$, $r>1$, be an Anosov volume preserving diffeomorphism
on a compact manifold $M$. Then $f$ is ergodic.
\end{thm}

\noindent {\it Proof}: Let $\Lambda \subset M$ be an invariant set and
$\mu(\Lambda) >0$. Let $\Lambda_s \subset \Lambda$ be the set such
that for each $x \in \Lambda_s$, a.e. ($\mu_s$) $y$ in $W^s_\delta(x)$
is in $\Lambda$. The above lemma shows that $\mu(\Lambda_s) = \mu
(\Lambda)$. Also by the above lemma, for a.e. $z \in \Lambda$,
a.e. ($\mu_u$) $x$ in $W^u_\delta(z)$ is in $\Lambda_s$. By the
absolute continuity of the stable and unstable foliations, the set
$$\bigcup_{y \in W^u_\delta(z)} W^s_\delta(y)$$ has the full measure in
the $\delta$ neighborhood of $z$ and therefore $\Lambda$ has the full
measure in a $\delta$ neighborhood of $z$ for a.e. $z \in
\Lambda$. This implies that $\Lambda$ has the full measure in
$M$. Since $\Lambda$ is an arbitrary positive measure set, this
implies ergodicity.

This proves the theorem.

Finally we remark that Lemma \ref{lemf} is also true for partially
hyperbolic invariant and non-uniformly hyperbolic invariant sets (sets with
non-zero Liapunov exponents). The proof is exactly the same.


\end{document}